\documentclass[preprint,authoryear,1p,times]{elsarticle}

\usepackage{amssymb}
\usepackage{amsthm}
\usepackage{amsmath,amsthm,amssymb,dsfont}
\usepackage{parallel,amsfonts,accents}
\usepackage{bm}
\newtheorem{thm}{Theorem}
\theoremstyle{plain}      \newtheorem{lem}{Lemma}
\theoremstyle{plain}      \newtheorem{dfn}{Definition}
\theoremstyle{plain}      
\theoremstyle{plain}      
\theoremstyle{definition} 
\theoremstyle{definition} 
\theoremstyle{definition} 
\theoremstyle{plain} \newtheorem{cor}{Corollary}
\theoremstyle{definition} 
\theoremstyle{plain} 
\theoremstyle{definition} 
\theoremstyle{definition} 
\theoremstyle{definition} 

\usepackage{pgf,tikz}
\usetikzlibrary{snakes}
\usetikzlibrary{backgrounds}
\usetikzlibrary{arrows}
\usetikzlibrary{patterns}
\usetikzlibrary{fit}
\usetikzlibrary{calc}

\newcommand\ds{\mathds}
\newcommand\bb{\mathbb}

\newcommand\arginf{\operatornamewithlimits{arginf}}

\newcommand\for{\,\,\forall\,\,}

\newcommand\lb{\lbrace}
\newcommand\lt{\left}

\newcommand\rb{\rbrace}
\newcommand\rt{\right}

\newcommand\tth{^\text{th}}

\newcommand\Na{\ds{N}} 
\newcommand\R{\bb{R}}  

\journal{Probability and Statistics Letters}
\begin{document}
\begin{frontmatter}



\title{Weighted Frechet Means as Convex Combinations in Metric Spaces:
           Properties and Generalized Median Inequalities.}

\author[label1,label2,label3]{Cedric E. Ginestet}
\author[label1,label2]{Andrew Simmons}
\author[label3]{Eric D. Kolaczyk}
\address[label1]{King's College London, Institute of Psychiatry,
  Department of Neuroimaging}
\address[label2]{National Institute of Health Research (NIHR) Biomedical
  Research Centre for Mental Health}
\address[label3]{Department of Mathematics and Statistics, Boston 
  University, MA}
\begin{abstract}
   In this short note, we study the properties of the weighted Frechet mean 
   as a convex combination operator on an arbitrary metric space
   $(\mathcal{Y},d)$. We show that this binary operator is commutative,
   non-associative, idempotent, invariant to
   multiplication by a constant weight and possesses an identity element.
   We also cover the properties of the weighted cumulative Frechet mean.
   These tools allow us to derive several types of median inequalities
   for abstract metric spaces that hold for both negative and positive Alexandrov
   spaces. In particular, we show through an example that these bounds
   cannot be improved upon in general metric spaces. For weighted
   Frechet means, however, such inequalities can solely be derived for weights
   equal or greater than one. This latter limitation highlights the
   inherent difficulties associated with abstract-valued random variables. 
\end{abstract}
\end{frontmatter}
\section{Introduction}\label{sec:intro}
The core task of statistics is to summarize data, which is commonly
done by identifying \textit{typical elements}. As observed by 
\citet{Frechet1948}, typical elements are the elements in the sample
space that are as similar as possible to all the other elements in that space. If a notion
of distance is defined on this space of elements, then it
follows that the most typical element is the one that differs the least
from all the others. Such an element is commonly referred to as the
Frechet mean, barycentre or Karcher mean \citep{Karcher1977}. These
typical elements have been studied in metric spaces at various
levels of generality. The almost sure convergence of the Frechet
sample mean to its theoretical analogue has been demonstrated for
separable metric spaces with bounded metric \citep{Ziezold1977}, in 
compact spaces \citep{Sverdrup1981}, and when the Frechet mean is
assumed to be unique \citep{Bhattacharya2003}. The concept of Frechet
mean has also proved to be useful in different domains of applications,
such as in image analysis \citep{Thorstensen2009,Bigot2011} or when
studying phylogenetic trees \citep{Balding2009}. 

One of the outstanding questions in this field is whether the
classical \textit{median inequality} can be
recovered using the Frechet mean. If a similar result can be derived,
it may then be possible to generalize standard results in probability
and statistics to abstract metric spaces. In the Euclidean plane, the 
median inequality states that for every triangle $\Delta ABC$, 
\begin{equation}
    EC \leq \frac{1}{2}\big(AC +  BC\big),
    \notag
\end{equation}
where $E$ is the midpoint of $\overline{AB}$. In this short note, 
we consider generalizations of this law to abstract metric spaces,
whereby the midpoint of a segment in such spaces is defined as the
Frechet mean of the endpoints of that segment. In particular, we
explore the properties of the Frechet mean and its cumulative
extension in arbitrary metric spaces for any type of
Alexandrov curvature \citep[see][for a study of typical elements in
negatively curved metric spaces]{Herer1992}. Our definition of the
Frechet mean as a convex combination bears some resemblance with the
convex combination operator introduced by \citet{Teran2006},
and we will draw some specific links with the work of these authors in
the sequel. 

\section{Weighted Frechet Mean}\label{sec:mean}
The Frechet mean generalizes the arithmetic mean to abstract metric spaces. In
general, this quantity may not be unique. When this is the case, the
set of all minimizers is referred to as the \textit{Frechet mean
  set}. Here, we will only consider a given element from
the set of all such minimizers. One of the interesting properties of the
Frechet mean is that it also allows the combination of
\textit{subsets} of a metric space $(\mathcal{Y},d)$. This constitutes another
important generalization of the classical notion of arithmetic
mean, where one solely combines elements of the Real line. We
therefore define the Frechet mean with respect to subsets of $(\mathcal{Y},d)$.
Here, the distance between a subset $A\subseteq\mathcal{Y}$ and a point
$y\in\mathcal{Y}$ is $d(y,A):=\inf\{d(y,a):a\in A\}$. 
\begin{dfn}\label{dfn:mean}
   On a given metric space $(\mathcal{Y},d)$, the Frechet
   mean of the $r\tth$ order is defined for any two subsets
   $A,B\subseteq\mathcal{Y}$ and real numbers $\alpha,\beta\geq 0$, as follows, 
   \begin{equation}
       \alpha A\oplus_{r} \beta B \in
       \arginf_{y\in\mathcal{Y}} \alpha d(A,y)^{r} + \beta d(y,B)^{r},
   \end{equation}
   for every $r\geq 1$. 
   Similarly, the Frechet cumulative mean operator of the $r\tth$
   order is defined for any finite sequence of subsets of $\mathcal{Y}$, denoted
   $A_{1},\ldots,A_{n}$, and sequence of non-negative real numbers
   $\alpha_{1},\ldots,\alpha_{n}$, as follows,
   \begin{equation}
        \sideset{}{_{r}}\bigoplus_{i=1}^{n}\alpha_{i}A_{i} \in
        \arginf_{y\in\mathcal{Y}} \sum_{i=1}^{n}\alpha_{i} d(A_{i},y)^{r}.
   \end{equation}
\end{dfn}
We now study the properties of the Frechet mean of any
order, in abstract metric spaces. The following lemma and corollary
are true for all $r\geq1$, and therefore this subscript is omitted. 
\begin{lem}\label{lem:mean}
   For any $A,B,C\subseteq\mathcal{Y}$ in a metric space $(\mathcal{Y},d)$,
   and for any $\alpha\geq 0$, the Frechet mean operator 
   satisfies the following:
   \begin{description}
   \item[(i.)] Commutativity: $A\oplus B = B\oplus A$;
   \item[(ii.)] Non-associativity: $A\oplus(B\oplus C) \neq (A\oplus B)\oplus C$; 
   \item[(iii.)] Idempotency: $A \oplus A=A$;
   \item[(iv.)]  Proportionality: $\alpha A \oplus \alpha B = A\oplus B$;
   \item[(v.)] Identity element: $A\oplus\mathcal{Y}=A$, for every $A$.
   \end{description}
\end{lem}
\begin{proof}
   Commutativity in (i) is immediate from the commutativity
   of addition on the real numbers. Invariance with respect to a constant weight in
   (iv) is a direct consequence of the definition of the Frechet
   mean operator. Idempotency in (iii) follows from
   the fact that metrics do not distinguish between identical
   elements, and therefore $d(A,y)^{r} + d(y,A)^{r}=2d(y,A)^{r}$, using the symmetry of
   $d$ and invoking (iv). The existence of the identity element in (v) can be deduced by
   noting that $d(A,y)+d(y,\mathcal{Y})$ is only null when $y\in A\cap\mathcal{Y}=A$.
   Finally, non-associativity in (ii) can be proved
   through a numerical counter-example in $\R$ equipped with the
   Euclidean metric (i.e.~ take any three distinct real numbers).
   Therefore, associativity does not hold in general. 
\end{proof}
Equivalent properties can be immediately deduced from lemma \ref{lem:mean} for
the case of the Frechet cumulative mean operator, as described in the
following corollary. 
\begin{cor}\label{cor:cumulative}
   For any sequence of subsets, denoted $\{A_{1},\ldots,A_{n}\}$, in a
   metric space $(\mathcal{Y},d)$, for any real number $\alpha\geq0$, 
   and any $n\in\Na$, the Frechet cumulative mean operator 
   satisfies for any label permutation $\nu:I\mapsto I$, where $I:=\{1,\ldots,n\}$, 
   \begin{description}
   \item[(i.)] Commutativity: $\bigoplus_{i=1}^{n} A_{i} = \bigoplus_{i=1}^{n}A_{\nu(i)}$;
   \item[(ii.)] Proportionality: $\bigoplus_{i=1}^{n}\alpha A_{i}=\bigoplus_{i=1}^{n}A_{i}$.
   \item[(iii.)] Idempotency: $\bigoplus_{i=1}^{n} A = A$;
   \end{description}
\end{cor}
Note that property (i) of the Frechet cumulative mean corresponds to
condition (i) in \citet{Teran2006}. These authors have studied the
behavior of convex combination operators in metric spaces. However,
the regrouping condition, denoted (ii) in \citet{Teran2006} does not
hold in general abstract spaces for the Frechet mean, due to its
non-associativity. This lack of associativity will also lead to some
difficulties when extending the generalized median inequality from the binary
Frechet mean to its cumulative analogue. 

\begin{dfn}\label{dfn:convex}
    A set $A\subseteq(\mathcal{Y},d)$ is $\oplus$-convex if for every sequence
    $y_{1},\ldots,y_{n}$ of points in $A$ and non-negative numbers
    $\alpha_{i}$'s, we have $\bigoplus_{i=1}^{n}\alpha_{i}y_{i} \in A$.
\end{dfn}
The Frechet mean operator therefore allows the construction of $\oplus$-\textit{convex hulls} in
$(\mathcal{Y},d)$, such that for every $A\subseteq \mathcal{Y}$, the $\oplus$-convex hull of
$A$ of the $r\tth$ order is defined as
\begin{equation}
     H^{r}(A) :=
    \lt\lb\lt. \sideset{}{_{r}}\bigoplus_{i=1}^{n}\alpha_{i}y_{i}\rt| 
    y_{i}\in A, \alpha_{i}\geq0, \for\!n\in\Na\rt\rb,
    \label{eq:hull}
\end{equation}
Here, although $H^{r}(A)$ is $\oplus$-convex, $A$ need not be convex
in the classical sense. That is, if $A$
is a subset of a vector space, for instance, it may not be convex with respect to
vector addition. Nonetheless, given any metric on that vector space, one can
construct a hull, which is convex with respect to the Frechet mean
based on that particular metric. 

By definition, $H^{r}(A)$ is $\oplus_{r}$-convex for every $r\geq1$.
Similarly, observe that the \textit{closure} of $\mathcal{Y}$ is trivially $\oplus_{r}$-convex. 
Although the definition in equation (\ref{eq:hull}) appears to be the
one of a convex cone, in fact, it defines a convex hull. 
That is, although in our adopted definition,
we have not explicitly required the $\alpha_{i}$'s to sum to
$1$, it follows from the \textit{proportionality} of the Frechet cumulative
mean that these weights can be normalized without altering the choice
of the optimal elements in $\mathcal{Y}$. Since by definition, the Frechet
cumulative mean of a collection of points, $\{y_{1},\ldots,y_{n}\}$,
is necessarily located in the convex hull of these points, it
follows that the Frechet cumulative mean of any order can be regarded
as a \textit{convex combination} on $(\mathcal{Y},d)$. Note, however, that this concept
is here used in a more general sense than in \citet{Teran2006}.

\section{Median Inequalities in Metric Spaces}\label{sec:median
  inequality}
We here state and prove the main results of this paper for the Frechet
mean and cumulative mean of the \textit{first} order. Hence,
in this section, all Frechet operations will be assumed to be
conducted with respect to $r=1$. The more general case will be studied
in section \ref{sec:r order}. Note also that, without loss of generality, we have
formulated these results in terms of single elements in
$(\mathcal{Y},d)$. However, all of these results also hold for subsets of $\mathcal{Y}$. 
\begin{thm}\label{thm:median inequality}
     For any abstract metric space $(\mathcal{Y},d)$, and for every $x,y,\xi\in\mathcal{Y}$,
     \begin{equation}
         d(\xi,x\oplus y) \leq d(x,\xi) + d(\xi,y).
         \notag
     \end{equation}
\end{thm}
\begin{proof}
     Assume that the result does not hold, and that for some
     $x,y,\xi\in\mathcal{Y}$, we have instead 
     $d(\xi,x\oplus y) > d(\xi, x) + d(\xi,y)$.
     By the triangle inequality with respect to $y$, it follows that 
     \begin{equation}
         d(x,\xi) + d(\xi,y) < d(\xi,x\oplus y) \leq d(\xi,y) +
         d(y,x\oplus y),
         \notag
     \end{equation}
     which simplifies to $d(x,\xi) < d(y,x\oplus y)$.
     Similarly, by invoking the triangle inequality with respect to
     $x$, we obtain $d(\xi, y) < d(x,x\oplus y)$.       
     Now, combining these two strict inequalities and using the
     symmetry of $d$, this gives
     \begin{equation}
         d(x,\xi) + d(\xi,y) < d(x,x\oplus y) + d(x\oplus y,y),
         \notag
     \end{equation}
     but this contradicts the minimality of $x\oplus y$, and therefore proves
     the theorem. 
\end{proof}
Observe that the Euclidean median law does not hold in general metric
spaces, as illustrated by figure \ref{fig:tight bound} for a
negatively curved Alexandrov space \citep[see][]{Burago2001}. Therefore,
the result in theorem \ref{thm:median inequality} is tight in the sense that
this inequality can be saturated for some metric spaces. By contrast,
the inequality is strict in the Euclidean case. A similar
inequality can be derived for the case of a weighted Frechet mean,
albeit observe that such weights should be equal or greater than 1.
\begin{figure}[t]
\centering
\tikzstyle{background rectangle}=[draw=gray!30,fill=gray!5,rounded corners=1ex]
\begin{tikzpicture}[font=\small,scale=.8]
    \draw(-1.5,0)node(xi){};
    \draw(2.8,3)node(x){};
    \draw(3,.75)node(xy){};
    \draw(3.7,-1.3)node(y){};
    \draw(xi)node[anchor=north east]{$\xi$};
    \draw(x)node[anchor=south west]{$x$};
    \draw(xy)node[anchor=west]{$x\oplus y$};
    \draw(y)node[anchor=north west]{$y$};
    \draw (xi) parabola bend (xi) (x);
    \draw (xi) parabola bend (xi) (y); 
    \draw (xy) parabola bend (xy) (x); 
    \draw (y) parabola bend (y) (xy); 
    \draw (xi) parabola bend (xi) (xy); 
    \fill[fill=red!50](x) circle (2.5pt); 
    \fill[fill=red!50](xy) circle (2.5pt); 
    \fill[fill=red!50](y) circle (2.5pt); 
    \fill[fill=red!50](xi) circle (2.5pt); 
    \draw(2,-.1)node[anchor=north]{$1$};
    \draw(1,1)node[anchor=south]{$1$};
    \draw(1.5,.3)node[anchor=south]{$2$};
    \draw(2.8,1.9)node[anchor=west]{$1$};
    \draw(3.2,-.4)node[anchor=west]{$1$};
\end{tikzpicture}
  \caption{Illustrative metric space with negative Alexandrov
    curvature \citep{Burago2001}. Albeit the distances between each point in this space
    satisfy the triangle inequality, we nonetheless have $d(\xi,
    x\oplus y) > \frac{1}{2}(d(\xi,x) + d(y,\xi))$, and therefore the
    classical Euclidean median inequality does not hold in this metric space. 
    In addition, observe that, in this setting, $\xi$ also constitutes a mean between
    $x$ and $y$, since we have $d(x,\xi)+d(\xi,y)=d(x,x\oplus y)+d(x\oplus
    y,y)$. 
   \label{fig:tight bound}}
\end{figure}
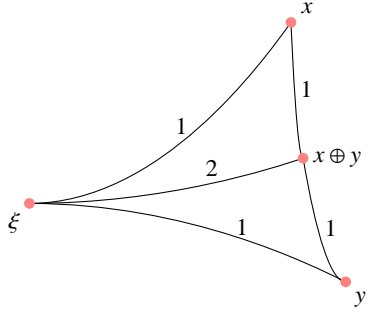  
\begin{cor}\label{cor:median inequality}
     For any abstract metric space $(\mathcal{Y},d)$, for every
     $x,y,\xi\in\mathcal{Y}$, and any $\alpha,\beta\geq1$,
     \begin{equation}
         d(\xi,\alpha x\oplus \beta y) \leq \alpha d(x,\xi) + \beta d(\xi,y).
         \notag
     \end{equation}
\end{cor}
\begin{proof}
    The proof is similar to the one of theorem \ref{thm:median
      inequality}, and also proceeds by contradiction. Assuming the reverse
    and invoking the triangle inequality and using the fact that
    $\alpha\geq1$, we have
    \begin{equation}
        \alpha d(x,\xi) + \beta d(\xi,y) < d(\xi,\alpha x\oplus \beta
        y) \leq \alpha d(\xi, x) + \alpha d(x,\alpha x\oplus \beta y),
          \notag
    \end{equation}
    which reduces to $\beta d(\xi,y)< \alpha d(x,\alpha x\oplus \beta
    y)$. Through an analogous procedure, we may obtain 
    $\alpha d(\xi,x)< \beta d(y,\alpha x\oplus \beta y)$, and
    combining these inequalities this gives the desired contradiction.
\end{proof}
Theorem \ref{thm:median inequality} can be generalized for the
Frechet cumulative mean operator. This result essentially states that
the cumulative mean operator is \textit{countably additive}. 
When extending these inequalities to the case of the cumulative
mean operator, observe that the non-associativity of the Frechet mean
does not allow a direct proof by induction, and therefore other arguments
have to be deployed in order to prove that an equivalent result holds
in this general setting. 
\begin{thm}[\textbf{Countable Additivity}]\label{thm:generalized median inequality}
     For any abstract metric space $(\mathcal{Y},d)$, for every sequence
     $\{y_{1},\ldots,y_{n}\}$ and $\xi$ in $\mathcal{Y}$, and for every $n\in\Na$,
     \begin{equation}
         d\lt(\bigoplus_{i=1}^{n}y_{i},\xi\rt) \leq 
         \sum_{i=1}^{n} d\lt(y_{i},\xi\rt).
         \notag
     \end{equation}
\end{thm}
\begin{proof}
     Again, seeking a contradiction, assume that
     $\sum_{i=1}^{n}d(y_{i},\xi) <
     d(\bigoplus_{i=1}^{n}y_{i},\xi)$. It then follows that
     through $n$ applications of the triangle inequalities, we obtain
     the following system of strict inequalities, 
     \begin{equation}
         \sum_{i=1}^{n}d(y_{i},\xi) <
         d\lt(\bigoplus_{i=1}^{n}y_{i},\xi\rt) \leq d(\xi,y_{k}) +
         d\lt(y_{k},\bigoplus_{i=1}^{n}y_{i}\rt),
       \notag
     \end{equation}
     for every $k=1,\ldots,n$. Each of these inequalities can be
     expanded using the positivity of $d$, as follows,
     \begin{equation}
         \sum_{i=1}^{n}d(y_{i},\xi) < d(\xi,y_{k}) +
         d\lt(y_{k},\bigoplus_{i=1}^{n}y_{i}\rt) 
         + \sum_{i\in I\setminus \{k,l\}} d(y_{i},\xi), 
       \notag
     \end{equation}
     where $I:=\{1,\ldots,n\}$, and where $l:= k+1$ if $k<n$ and
     $l:=1$ if $k=n$. Since the latter inequality holds for
     every $k=1,\ldots,n$, it suffices to sum these $n$ inequalities
     in order to obtain 
     \begin{equation}
         n\sum_{i=1}^{n}d(y_{i},\xi) - \sum_{k=1}^{n}d(\xi,y_{k})  
         - \sum_{k=1}^{n}\sum_{i\in I\setminus \{k,l\}} d(y_{i},\xi)
          < \sum_{k=1}^{n} d\lt(y_{k},\bigoplus_{i=1}^{n}y_{i}\rt),
         \notag
     \end{equation}
     which leads to $\sum_{i=1}^{n}d(y_{i},\xi) < \sum_{k=1}^{n}
     d(y_{k},\bigoplus_{i=1}^{k} y_{i})$.
     However, this contradicts the minimality of $\bigoplus_{i=1}^{k}
     y_{i}$, as desired. 
\end{proof}
\begin{cor}\label{cor:generalized median inequality}
     For any abstract metric space $(\mathcal{Y},d)$, for every sequence
     $\{y_{1},\ldots,y_{n}\}$ and $\xi$ in $\mathcal{Y}$, and for every
     sequence of real numbers $\{\alpha_{1},\ldots,\alpha_{n}\}$,
     satisfying $\alpha_{i}\geq1$, 
     \begin{equation}
         d\lt(\bigoplus_{i=1}^{n}\alpha_{i}y_{i},\xi\rt) \leq 
         \sum_{i=1}^{n} \alpha_{i}d\lt(y_{i},\xi\rt),
         \notag
     \end{equation}
     for every $n\in\Na$.
\end{cor}
\begin{proof}
     The proof strategy is similar to the one of theorem \ref{thm:generalized
       median inequality}, but using the argument described in the
     proof of corollary \ref{cor:median inequality}. That is, using
     the same notation as in the proof of theorem \ref{thm:generalized
       median inequality}, since $\alpha_{i}\geq1$ for all $\alpha_{i}$, one can derive the
     following system of $n$ inequalities, 
     \begin{equation}
         \sum_{i=1}^{n}\alpha_{i}d(y_{i},\xi) < \alpha_{k}d(\xi,y_{k}) +
         \alpha_{k}d\lt(y_{k},\bigoplus_{i=1}^{n}\alpha_{i}y_{i}\rt) 
         + \sum_{i\in I\setminus \{k,l\}} \alpha_{i}d(y_{i},\xi),
         \label{eq:proof1}
     \end{equation}
     Combining these inequalities gives $\sum_{i=1}^{n}\alpha_{i}d(y_{i},\xi) < \sum_{k=1}^{n}
     \alpha_{k}d(y_{k},\bigoplus_{i=1}^{k} \alpha_{i}y_{i})$,
     which provides the required contradiction.
\end{proof}

\section{Median Inequalities of the $r\tth$ Order}\label{sec:r order}
More generally, one may be interested in considering whether analogues
of the above median inequalities also hold for Frechet means and cumulative
means of arbitrary orders, i.e. for which $r\geq1$. The following two
results state such generalized versions of the median inequality. 
\begin{thm}\label{thm:r median inequality}
     For any abstract metric space $(\mathcal{Y},d)$, for every
     $x,y,\xi\in\mathcal{Y}$, and for every $r\geq1$, the Frechet mean of
     the $r\tth$ order satisfies,
     \begin{equation}
         d(\xi,x\oplus_{r} y)^{r} \leq 2^{r-1}\Big(d(x,\xi)^{r} + d(\xi,y)^{r}\Big).
         \notag
     \end{equation}
\end{thm}
\begin{proof}
     Using the same argument described in the proof of theorem
     \ref{thm:median inequality}, we assume for contradiction that 
     $d(\xi,x\oplus_{r} y)^{r} > 2^{r-1}(d(x,\xi)^{r} +
     d(\xi,y)^{r})$ holds. Here, we will require a result due to
     \citet{Frechet1948}, which states that 
     \begin{equation}
           d(a,b)^{r} \leq 2^{r-1}\Big(d(a,c)^{r} + d(c,b)^{r}\Big),
           \notag
     \end{equation}
     for every $a,b,c\in\mathcal{Y}$ and every $r\geq1$. See equation (5) on page 228 of
     \citet{Frechet1948}. This equation will be referred to in the
     sequel as the \textit{triangle inequality of the $r\tth$ order}. By
     using this result, it immediately follows that 
     \begin{equation}
          2^{r-1}\Big(d(x,\xi)^{r} + d(\xi,y)^{r}\Big) 
          < d(\xi,x\oplus_{r} y)^{r} 
          \leq 2^{r-1}\Big(d(\xi,x)^{r} + d(x,x\oplus_{r} y)^{r}\Big),
          \notag
     \end{equation}
     which reduces to $d(\xi,y)^{r} <  d(x,x\oplus_{r}
     y)^{r}$. Similarly, we have $d(x,\xi)^{r} < d(y,x\oplus_{r}
     y)^{r}$. As before, combining these results 
     contradicts the minimality of $x\oplus_{r} y$, and therefore the
     result is true for the Frechet mean of the $r\tth$ order. But
     $r$ was arbitrary and thus the theorem holds for any $r\geq1$.
\end{proof}
It is straightforward to generalize this result to the weighted Frechet
cumulative mean of the $r\tth$ order. In its most general form, we
therefore have the following median inequality. 
\begin{cor}\label{cor:r generalized median inequality}
     For any abstract metric space $(\mathcal{Y},d)$, for every sequence
     $\{y_{1},\ldots,y_{n}\}$ of elements in $\mathcal{Y}$, and for every
     sequence of real numbers $\{\alpha_{1},\ldots,\alpha_{n}\}$,
     satisfying $\alpha_{i}\geq1$, 
     \begin{equation}
         d\lt(\sideset{}{_{r}}\bigoplus_{i=1}^{n}\alpha_{i}y_{i},\xi\rt)^{r} \leq 
         2^{r-1}\sum_{i=1}^{n} \alpha_{i}d\lt(y_{i},\xi\rt)^{r},
         \notag
     \end{equation}
     for every $n\in\Na$, for every $r\geq1$.
\end{cor}
\begin{proof}
     Using the arguments invoked in the proofs
     of the aforementioned results, we proceed by contradiction and
     assume that the reverse of the conclusion of corollary
     \ref{cor:r generalized median inequality} holds. We have by the
     triangle inequality of the $r\tth$ order, for every $k=1,\ldots,n$,
     \begin{equation}
         2^{r-1}\sum_{i=1}^{n}\alpha_{i}d(y_{i},\xi)^{r} <
         d\lt(\sideset{}{_{r}}\bigoplus_{i=1}^{n}\alpha_{i}y_{i},\xi\rt)^{r}
         \leq 
         2^{r-1}\lt[
         d\lt(\xi,y_{k}\rt)^{r} + 
         d\lt(y_{k},\sideset{}{_{r}}\bigoplus_{i=1}^{n}\alpha_{i}y_{i}\rt)^{r}
         \,\rt].
       \notag         
     \end{equation}     
     This can be expanded using the fact that $\alpha_{i}\geq1$, for
     every $i=1,\ldots,n$, and simplified by dividing both sides by
     $2^{r-1}$ in order to obtain the analogue of equation
     (\ref{eq:proof1}) but where all metrics are elevated to the
     $r\tth$ power and the Frechet cumulative mean is of the $r\tth$ order.
     As in the proof of theorem
     \ref{thm:generalized median inequality}, combining this system of
     $n$ strict inequalities contradicts the minimality of the
     Frechet cumulative mean of the $r\tth$ order, and this  completes
     the proof. 
\end{proof}

\bibliographystyle{elsarticle-harv}
\bibliography{/home/cgineste/ref/bibtex/Statistics,%
              /home/cgineste/ref/bibtex/Neuroscience}
\end{document}